  \documentclass[11pt]{article}
\begin{document}
\renewcommand{\theequation}{\arabic{section}.\arabic{equation}}
   \newtheorem{lemma}{Lemma}[section]
   \newtheorem{theorem}[lemma]{Theorem}
   \newtheorem{remark}[lemma]{Remark}
   \newtheorem{prop}[lemma]{Proposition}
   \newtheorem{coro}[lemma]{Corollary}
   \newtheorem{definition}[lemma]{Definition}
   \newcommand{\eps}{\varepsilon}
   \newcommand{\Wsob}{\smash{{\stackrel{\circ}{W}}}_2^1(D)}
   \newcommand{\EX}{{\mathbb E }}
   \newcommand{\PX}{{\mathbb P }}
\newcommand{\be}{\begin{equation}}
   \newcommand{\ee}{\end{equation}}
   \newcommand{\tr}{\triangle}
\newcommand{\e}{\epsilon}
\renewcommand{\a}{\alpha}
\renewcommand{\b}{\beta}
\newcommand{\om}{\omega}
\newcommand{\Om}{\Omega}
\newcommand{\D}{\Delta}
\newcommand{\p}{\partial}
\newcommand{\de}{\delta}
\renewcommand{\phi}{\varphi}
\newcommand{\N}{{\mathbb N}}
\newcommand{\R}{{\mathbb R}}
\newcommand{\cF}{{\cal F}}
\def\rit{I\!\!R}

\baselineskip = 16pt

\parindent0.0em
%%%%%%%%%%%%%%%%%%%%%%%%%%% %%%%

\title
{Dynamics of a Coupled Atmosphere-Ocean Model}

\author{Hongjun Gao$^1$ and Jinqiao Duan$^2$  \\ \\
1. Department of Mathematics\\Nanjing Normal University
\\Nanjing 210097,  China \\
and \\
Institute for Mathematics and its Applications\\
    University of Minnesota\\
    Minneapolis, MN 55455, USA
 \\ \\
2. Department of Applied Mathematics\\ Illinois Institute of Technology \\
  Chicago, IL 60616, USA \\E-mail:  duan@iit.edu \\
  and \\
Department of Mathematics\\University of Science and Technology of
China
\\Hefei 230026,  China
   }

\date{June 27, 2003/Revised December 29, 2003}

\maketitle

\begin{abstract}
{\bf Nonlinear Analysis B: 2004, in press}.

 We consider a coupled
atmosphere-ocean model, which involves hydrodynamics,
thermodynamics and nonautonomous interaction at the air-sea
interface. First, we show that the coupled atmosphere-ocean system
is stable under the external fluctuation in the atmospheric energy
balance relation. Then, we estimate the atmospheric temperature
feedback in terms of the freshwater flux, heat flux and the
external fluctuation at the air-sea interface, as well as the
earth's longwave radiation coefficient and the shortwave solar
radiation profile. Finally, we prove that the coupled
atmosphere-ocean system has time-periodic, quasiperiodic and
almost periodic motions, whenever the external fluctuation  in the
atmospheric energy balance relation is time-periodic,
quasiperiodic and almost periodic, respectively.

\vspace{0.9cm}


{\bf Abbreviated Title}:  Coupled Atmosphere-Ocean Dynamics

\end{abstract}

\newpage

%%%%      Section 1
%%%%%%%%%%%%%%%%%%%%%%%%%%%%%%%%%%
\section{Introduction: A coupled atmosphere-ocean model}

The global ocean circulation consists of the wind-driven upper ocean
circulation and
  a  meridional overturning deep ocean circulation called the
thermohaline circulation.  The ocean thermohaline
circulation involves water masses sinking at high latitudes and
upwelling at lower latitudes. During the thermohaline circulation, water masses carry heat  or cold
around the globe. Thus, it is believed that the global ocean thermohaline
circulation plays an important role in the climate \cite{Siedler}.

The thermohaline circulation   is maintained by water
density contrasts in the ocean, which themselves are created by
atmospheric forcing, namely, heat and freshwater exchange via
evaporation and  precipitation at the air-sea interface.
Thus the thermohaline circulation is
described by coupled  atmosphere-ocean  models \cite{PeiOor92, Siedler}.
Such coupled models also
describe  feedback of the thermohaline circulation on the atmospheric
dynamics (e.g., temperature feedback).

Mathematical models
 are a key component of our
understanding of  climate and geophysical  systems.
The
formulation and analysis of mathematical models is central to the
progress of better understanding of the  thermohaline circulation
dynamics and its impact on climate change.

We   consider a coupled  atmosphere-ocean model, with simplified
atmospheric dynamics, i.e., the atmospheric dynamics is described
by an energy balance model.

This is a zonally averaged, coupled atmosphere-ocean model on the
meridional, latitude-depth $(y,z)$-plane as used by various
authors \cite{Stocker, Wright, Chen, Dijk, DijkBook, DuanGaoSchm}.
This model has been shown to capture some interesting climate
phenomena \cite{Stocker, Wright, Chen}. It is composed of a
one-dimensional
 stochastic energy
balance model   proposed by North and Cahalan \cite{North},
for the    latitudinal
 atmosphere surface temperature  $\theta(y,t)$, together with
the Boussinesq equations for ocean dynamics in terms of stream
function $\psi(y,z,t)$,
and transport  equations for the oceanic salinity $S(y, z, t)$
and the oceanic temperature $T(y, z, t)$ on
the domain $D=\{ (y,z): 0\leq y, z \leq 1\}$:
\begin{eqnarray}
\theta_t   =   \theta_{yy} -(a+\theta) &+& S_a(y)-\gamma(y)[S_o(y)+\theta -T]
        + f(y, t),  \:0\le y \le 1,\label{eqn1}\\
q_t + J(q, \psi )      &  = & Pr \D q  +PrRa (T_y - S_y), \;\;\; (y, z) \in D, \label{eqn2}\\
T_t + J(T, \psi )       &  = & \D T, \;\;\;\;     (y, z) \in D,\label{eqn3}\\
S_t + J(S, \psi )        &  = & \D S,  \;\;\;\;     (y, z) \in D,\label{eqn4}
\end{eqnarray}
where   $q(y, z, t)= -\D \psi$ is the vorticity;
  velocity field is $(v, w)= (\psi_z, -\psi_y)$;  $a$ is a   positive constant
parameterizing the effect
of the earth's longwave radiative cooling; $S_a(y)$ and $ S_o(y)$
are empirical functions representing the latitudinal dependence of the shortwave solar radiation;
$\gamma(y)$  is  the latitudinal  fraction of the earth covered
by the ocean basin;  Pr is the
Prandtl number  and Ra is the Rayleigh number.  The first equation is the
energy balance model proposed by
North and Cahalan \cite{North}. The   forcing
$f(y, t)$ may arise from, for example, eddy transport fluctuation, stormy
bursts of latent heat, and flickering
 cloudiness variables.  Moreover, $J(g,h)=g_xh_y-g_yh_x$
is the Jacobian operator and $\D=\p_{yy}+\p_{zz}$ is the Laplace
operator.  The effect of the rotation   is
parameterized in the magnitude of the viscosity and diffusivity terms
as discussed in \cite{Thual}.

The  no-flux boundary condition  is taken for
the  atmosphere  temperature $\theta (y, t)$
\begin{equation}\label{btheta}
\theta_y (0, t)=\theta_y (1, t)  = 0.
\end{equation}
The  fluid boundary condition is no normal flow and free-slip
on the whole boundary
\begin{equation}\label{bvorticity}
\psi=0,  \;  q = 0.
\end{equation}

The flux boundary conditions are assumed
for the ocean temperature $T$ and salinity $S$.

 At top $z=1$, the fluxes are specified as:
\begin{equation}\label{btop}
 T_z = S_o(y)+ (\theta-T)|_{z=1}, \; \;    S_z=F(y),
 \end{equation}
with $F(y)$ being the given freshwater flux.

At bottom $z=0$:
\begin{equation}\label{bbottom}
T_z =  S_z=0  \;.
\end{equation}

On the lateral boundary    $y=0, 1$:
\begin{equation}\label{side}
T_y =S_y =0.
\end{equation}
We also assume the following compatibility condition:
\begin{equation}\label{comp}
S_o^{\prime}(0) = S_o^{\prime}(1) = F^{\prime}(0) = F'(1) = 0.
\end{equation}

The non-autonomous partial differential equation (\ref{eqn1})
is only defined on the air-sea
interface and it may be regarded as a   dynamical
boundary condition.  The boundary condition (\ref{btop}) involves
a coupling between the atmospheric and oceanic temperature at the
air-sea interface.

In the  next section, we discuss the well-posedness
of this coupled atmosphere-ocean model. Then we investigate
the stability  of this coupled system
under external fluctuation in \S 3, atmospheric temperature feedback
in \S 4, time-periodic, quasiperiodic and almost periodic coupled motion
in \S 5, respectively. Finally, we summarize these results  in the final section
\S 6.

%%%%%%%%%%%% Section 2
%%%%%%%%%%%%%
\section{Mathematical Setup}

In order to use the standard result in \cite{henry} for the local existence, we
 homogenize inhomogeneous
 boundary conditions
for $T, S$ on the top boundary $z=1$ as in \cite{Lions92}.

First, we construct two scalar functions:
$$T_{\e}^* = \tilde{T}^*\eta_{\e}(z), \;\;S_{\e}^* = \tilde{S}^*\eta_{\e}(z),
\;\forall \e\in(0, \frac12),$$
where
\begin{eqnarray*}
\tilde{T}^* &=& [S_o(y) + \theta](1 - e^{1 - z}),
\\
\tilde{S}^* &=& F(y)z,\\
\eta_{\e}(z) &\in& C^{\infty}([0, 1])\;\mbox{is given by}
\end{eqnarray*}

\[ \eta_{\e}(z) = \left\{\begin{array}{cc}
 1, & 1 - \e \le z \le 1,\\
\mbox{increasing}, & 1 - 2\e \le z \le 1 - \e,\\
0,  & 0 \le z \le 1 - 2\e.
\end{array}
\right.
\]
Then, set
$$\hat{T} = T - T_{\e}^*, \;\; \hat{S} = S - S_{\e}^*.
$$
By (\ref{btheta}) and (\ref{comp}), we see that the boundary conditions (\ref{btop})
for the new variables  $\hat{T}$ and $\hat{S}$ are homogenized and do not
affect other boundary conditions.
Thus  (\ref{eqn1})--(\ref{side}) become (for the simplicity, we
still use
$T$ and $S$
instead of $\hat{T}$ and $\hat{S}$)

\be
\theta_t   =   \theta_{yy} -(a + \theta) + S_a(y) - \gamma(y)[S_o(y)+\theta -T]
        + f(y, t), \;  0\le y \le 1,\label{neqn1}
\ee
\be
q_t + J(q, \psi )      =  Pr \D q  + Pr Ra (T_y - S_y + T_{\e y}^* - S_{\e y}^*),
\;\;\;  (y, z) \in D, \label{neqn2}
\ee
$$
T_t + J(T, \psi )  +   J(T_{\e}^*, \psi )    =  \D T +$$
$$
[((1 - e^{1-z})(1 - \gamma(y)) - e^{1-z})\eta_{\e}(z) + (1 - e^{1-z})\eta''_{\e}(z) + 2e^{1-z}\eta'_{\e}(z) - e^{1-z}\eta_{\e}]\theta$$
\be -
 (1 - e^{1-z})\eta_{\e}(z)\gamma(y)T(y, 1) + g, \;\;\;  (y, z) \in D,\label{neqn3}
 \ee
 $$
S_t + J(S, \psi )   +   J(S_{\e}^*, \psi )      = \D S $$
\be
+ F''(y)z\eta_{\e}(z) + F(y)(2\eta'_{\e}(z) + z \eta_{\e}''(z)), \;\;\;  (y, z) \in D,\label{neqn4}
\ee
where
$$
g(y, z, t) = - (1 - e^{1-z})\eta_{\e}(z)\{-a + S_a(y) - \gamma(y)S_o(y) + f(y, t) - S''_o(y)\}$$
$$ + [(1 - e^{1-z})\eta''_{\e}(z) + 2e^{1-z}\eta'_{\e}(z) - e^{1-z}\eta_{\e}]S_o(y).$$

The corresponding boundary conditions become:

On the whole boundary, the fluid flow satisfies
\begin{equation}\label{bvorticity1}
\psi=0,  \; \; q = 0.
\end{equation}

The boundary conditions for the atmosphere  temperature $\theta (y, t)$
(defined only on the air-sea interface) are
\begin{equation}\label{btheta1}
\theta_y (0, t)=\theta_y (1, t)  = 0.
\end{equation}

The boundary conditions for the ocean temperature $T$ and salinity $S$ become

  At top $z=1$:
\begin{equation}\label{btop1}
 T_z + T|_{z=1} = 0; \;\;   S_z = 0.
 \end{equation}

At bottom $z=0$:
\begin{equation}\label{bbottom1}
T_z =  S_z=0  \;.
\end{equation}

At  the lateral boundary $y=0, 1$:
\begin{equation}\label{side1}
T_y =S_y =0.
\end{equation}
The appropriate initial data $\theta_0, q_0, T_0, S_0$  are  also assumed.

Using the theory in \cite{henry}, we can
 obtain the following local existence theorem for problem (2.1)--(2.9) (that is (1.1)--(1.9)).

 \begin{theorem}({\bf Local Well-Posedness})
 \label{local}
 Let $\theta_0 \in H^1(0, 1)$,
$q_0\in H_0^1(D)$, $T_0, S_0\in H^1(D)$, $f\in L^{\infty}(0, \infty; L^2(0,
 1))$. Assume  that the physical data satisfy
$\gamma(y) \in L^{\infty}(0, 1)$, and $S_o(y), F(y) \in H^2(0, 1)$
and aslo assume that the compatibility condition (\ref{comp}) be satisfied.
 Then the coupled atmosphere-ocean system
 (\ref{neqn1})--(\ref{side1}) (that is (\ref{eqn1})--(\ref{side}) )
has a unique  (The uniqueness of $S$ is up to a constant)
local solution satisfying
 $$\theta\in L^{\infty}(0, \tau; H^1(0, 1))\cap L^2(0, \tau; H^2(0, 1)),$$
 $$
 q \in L^{\infty}(0, \tau; H_0^1(D))\cap L^2(0, \tau; H^2(D)\cap H_0^1(D)),$$
 $$T, S \in L^{\infty}(0, \tau; H^1(D)\times H^1(D))\cap L^2(0, \tau; H^2(D)\times
 H^2(D)),$$
 where $\tau$ depends on initial data and physical data $S_a(y), S_o(y), F(y)$ and
$f(y,t)$.
 \end{theorem}

 Since $-\D \psi =  q \in H_0^1$, we get $\psi\in H_0^1(D)\cap H^3(D)$. Hence
 the Jacobian  $J(\cdot, \cdot)$ is continuous from $H^1(D)\times H^3(D)\to
 L^2(D)\times L^2(D)$.

 \vspace{0.5cm}
 In order to obtain the global existence, we need a priori estimates. First,  we give a priori estimates for
 (\ref{eqn1})--(\ref{side}) in $L^2$.  In the sequel, $\|\cdot\|$ and $\|\cdot\|_1$ denote the norm of $L^2 $ and $H^1$  respectively.

Multiplying (\ref{eqn1}) by $\theta$ and performing the integration by parts,
  we conclude that
$$\frac12\frac{d}{dt}\|\theta\|^2 = - \|\theta_y\|^2 -
\|\theta\|^2 $$

\be
- a\int_0^1\theta dy + \int_0^1S_a\theta dy -
\int_0^1r(y)[S_o + \theta - T]\theta dy + \int_0^1f\theta dy,
 \ee

By the Cauchy-Schwarz inequality, we arrive at
\be \label{thetal2}
\frac12\frac{d}{dt}\|\theta\|^2 \le  - \|\theta_y\|^2 + C_1\|\theta\|^2  + \e_1\|T(y, 1)\|^2 + M_1,
\ee
where constant $M_1$   depends on $\|S_a||, \|S_o\|, a$ and
$\sup\limits_{0 \le t < \infty}\|f\|$, constant  $C_1$ depends on $\|\gamma\|_{L^{\infty}}$
and $\e_1>0$, $\e_1>0$ will be chosen later.

Multiplying (\ref{eqn2}) by $q$, performing the integration by parts and
 using the property of Jacobian and (\ref{bvorticity}), we have
 \be\label{ql2}
\frac12\frac{d}{dt}\|q\|^2 = - Pr\|\nabla q\|^2 +
PrRa\int_{D}(T_y - S_y)q.
 \ee

Similarly, from (\ref{neqn3}) and  (\ref{neqn4}),  we have
\be\label{tl2}
\frac12\frac{d}{dt}\|T\|^2 = - \|\nabla T\|^2 + \int_0^1[S_o +
\theta - T(y, 1)]T(y, 1) dy , \ee
and
 \be\label{sl2}
\frac12\frac{d}{dt}\|S\|^2 = - \|\nabla S\|^2 + \int_0^1F(y)S(y,
1) dy.
 \ee
 Note that
 $$
 PrRa\int_{\Omega}(T_y - S_y)q \le \frac{Pr}{2\lambda_1}\|q\|^2 +
 \frac{PrRa^2\lambda_1}{2}(\|T_y\|^2 + \|S_y\|^2)$$
 $$
  \le \frac{Pr}{2}\|\nabla q\|^2 +
 PrRa^2\lambda_1(\|\nabla T\|^2 + \|\nabla S\|^2),
 $$
 where $\lambda_1$ is a constant in the inequality
$\|v\|^2 \le \lambda_1\| \nabla v\|, v\in H_0^1$. Thus (\ref{ql2}) can be rewritten as
 \be\label{ql2ineq}
 \frac12\frac{d}{dt}\|q\|^2 \le - \frac{Pr}{2}\|\nabla q\|^2  +
 {PrRa^2\lambda_1}(\|\nabla T\|^2 + \|\nabla S\|^2).
 \ee
 Multiplying    (\ref{tl2}) by  $2PrRa^2\lambda_1$ and (\ref{sl2}) by
 $2PrRa^2\lambda_1$, and adding to  (\ref{ql2ineq}), we get
 $$
 \frac12\frac{d}{dt}(\|q\|^2 + 2PrRa^2\lambda_1(\|T\|^2 + \|S\|^2))
   $$
   $$
   \le - \frac{Pr}{2}\|\nabla q\|^2  -
 PrRa^2\lambda_1(\|\nabla T\|^2 + \|\nabla S\|^2)$$
 \be\label{qtsineq}
  + 2PrRa^2\lambda_1\int_0^1[S_o +
\theta - T(y, 1)]T(y, 1) dy +
2PrRa^2\lambda_1\int_0^1F(y)S(y, 1) dy.
 \ee
By the Cauchy-Schwarz inequality and  the trace inequality(\cite{evans}), we have
$$
 \frac12\frac{d}{dt}(\|q\|^2 + 2PrRa^2\lambda_1\|T\|^2 + 2PrRa^2\lambda_1\|S\|^2)$$
 $$
   \le - \frac{Pr}{2}\|\nabla q\|^2  - PrRa^2\lambda_1(\|\nabla T\|^2 + \|\nabla S\|^2)$$
 \be\label{qtsienq1}
  + C_2\|\theta\|^2 - \frac12\|T(y,1)\|^2 + \e_2(\|\nabla S\|^2 + \|S\|^2) + M_2,
 \ee
 where $C_2$ depends on $Pr, Ra$ and $\lambda_1$,  and
$M_2$ depends on $Pr, Ra, \lambda_1, \|S_o\|$ and $\|F\|$.
 Choosing $\e_1 < \frac12$ and $\e_2 < \frac{PrRa^2\lambda_1}{2}$,
combining (\ref{thetal2}) with (\ref{qtsienq1}), we obtain
 $$
 \frac12\frac{d}{dt}(\|\theta\|^2 + \|q\|^2 + 2PrRa^2\lambda_1\|T\|^2 + 2PrRa^2\lambda_1\|S\|^2) $$
 \be\label{totalineq}
  \le - \alpha(\|\nabla\theta\|^2 + \|\nabla q\|^2 + \|\nabla T\|^2 + \|\nabla S\|^2) + C_3(\|\theta\|^2 + \|S\|^2) + M_3,
 \ee
 where $C_3$ depends on $C_1$ and $C_2$,
$M_3$ depends on $M_1$ and $M_2$,  and
$\alpha$ is a positive constant   depending on $Pr, Ra$ and $\lambda_1$.
 By  the Gronwall inequality, we have
 $$
 \|\theta\|^2 + \|q\|^2 + \|T\|^2 + \|S\|^2 $$
 \be\label{l2estimate}
  + \int_0^b(\|\nabla\theta\|^2 + \|\nabla q\|^2 + \|\nabla T\|^2 + \|\nabla S\|^2)dt \le C_1(b),
 \ee
 for any given future
 time $b (0 < b < \infty, 0 < t \le b)$ and some positive constant $C_1(b)$ depending
 on $b, C_3$ and $M_3$. By a similar argument in \cite{gaoduan} (here we need to obtain
 the estimates in $H^1\times H^1 \times H^1 \times H^1 \times H^1$ for system of (\ref{neqn1})-(\ref{neqn4})  in order to avoid the trouble of
 non-homogeneous boundary conditions,
 we omit the details here since we will give a similar proof in
 \S 4), we have
 $$
 \|\nabla\theta\|^2 + \|\nabla q\|^2 + \|\nabla T\|^2 + \|\nabla S\|^2 + $$
 \be\label{h1estimate}
 \int_0^b(\|\D\theta\|^2 + \|\D q\|^2 + \|\D T\|^2 + \|\D S\|^2)dt \le C_2(b),
 \ee
 for any given $b (0 < b < \infty)$ and some positive constant $C_2(b)$
  depending on $\|S_o\|_{H^2}$, $\|F\|_{H^2}$, $b$ and $C_1(b)$.

 \begin{remark}
 In fact,  the estimate  we get  in (\ref{h1estimate})  is
 for $\hat{T}$ and $\hat{S}$, from which we can get the estimate
 for oringinal $T$ and $S$ using the estimate for $\theta$.
 \end{remark}

 \begin{remark}
Since the
 equivalence $\|T\|^2 + \|\D T\|^2$ with $\|T\|_{H^2}^2$  and the same for $S$, together with (\ref{l2estimate}),
 we can replace $\int_0^b(\|\D\theta\|^2 + \|\D q\|^2 + \|\D T\|^2 + \|\D S\|^2)dt$
 by $\int_0^b(\|\theta\|_{H^2}^2 + \|q\|_{H^2}^2 + \|T\|_{H^2}^2 + \|S\|_{H^2}^2)dt$
in the estimate  (\ref{h1estimate}).
 \end{remark}

 \vspace{0.3cm}
 With these global estimates, we have the following global existence theorem
for the coupled atmosphere-ocean system:

 \begin{theorem}
 \label{global}({\bf Global  Well-Posedness})
 Let $\theta_0 \in H^1(0, 1)$,
$q_0\in H_0^1(D)$, $T_0, S_0\in H^1(D)$, $f\in L^{\infty}(0, \infty; L^2(0,
 1))$. Assume  that the physical data satisfy
$\gamma(y) \in L^{\infty}(0, 1)$, and $S_o(y), F(y) \in H^2(0, 1)$
and aslo assume that the compatibility condition (\ref{comp}) be satisfied.
 Then for any given $b  (0 < b < \infty)$,
  the  coupled atmosphere-ocean system
 (\ref{neqn1})--(\ref{side1}) (that is (\ref{eqn1})--(\ref{side})) has a unique
(The uniqueness of $S$ is up to a constant) global solution satisfying
 $$\theta\in L^{\infty}(0, b; H^1(0, 1))\cap L^2(0, b; H^2(0, 1)),$$
 $$
 q \in L^{\infty}(0, b; H_0^1(D))\cap L^2(0, b; H^2(D)\cap H_0^1(D)),
 $$
  $$
  T, S \in L^{\infty}(0, b; H^1(D)\times H^1(D))\cap L^2(0, b; H^2(D)\times
 H^2(D)).
 $$
 \end{theorem}

 In the rest of this paper, we assume the conditions in this theorem
 are satisfied, so that we always have global unique solutions.

In the next section,  we  consider the stability of the above
coupled atmosphere-ocean system with respect to the
external fluctuation  $f(y, t)$  in the atmospheric energy balance
dynamics (\ref{eqn1}).

%%%%%%%%%%%%%%%%%%%%%%%%%%%%% section 3
%%%%%%%%%%%%%%%%%%%%%%%%%%%%%
\section{Stability under External Fluctuation}

\setcounter{equation}{0}

Paleo-evidence on the instability of    the thermohaline circulation
is now abundant.
Numerical work suggested that a sufficiently large external forcing
(such as external fluctuations in the atmospheric energy balance model
and the freshwater flux at the air-sea interface)
could destabilize   or shutdown  the thermohaline circulation
\cite{Rahms2000}.  This indicates that current
capacity of  carrying heat poleward by the  thermohaline circulation
may change when the  freshwater budget is altered.
Since the  thermohaline circulation's important role in redistributing
the heat around the globe,
a breakdown or instability of the current  thermohaline circulation
may lead to dramatic climate change \cite{Siedler}.
Because of this  relation between  the
thermohaline circulation and climate change,
there is growing interest in its stability or instability.
This   makes  the
stability issue of  the  thermohaline circulation not only of
scientific but also of great practical importance.

In this section, we prove the stability of the coupled atmosphere-ocean
system with respect to the
external fluctuation  $f(y, t)$  in the atmospheric energy balance
dynamics (\ref{eqn1}), i.e., the continuous dependence of solution
on $f(y, t)$ in
the space $H^1$.

\bigskip

Assume that $\{\theta_1, q_1, T_1, S_1\}$ and
$\{\theta_2, q_2, T_2, S_2\}$ are solutions
with respect to $f_1(y,t)$ and $f_2(y,t)$. Let
$$\bar \theta = \theta_1 - \theta_2, \; \bar q = q_1 - q_2, \; \bar T = T_1 - T_2, \; \bar S = S_1 - S_2, \; \bar f = f_1 - f_2,$$
then $\bar \theta, \bar q, \bar T$ and $\bar S$ satisfy
\begin{eqnarray}
\bar\theta_t   =   \bar\theta_{yy} - \bar\theta &-&\gamma(y)[\bar\theta - \bar T] + \bar f(y, t), \;\;\;\;\;\;\;\;\; 0\le z \le 1,\label{nneqn1}\\
\bar q_t + J(q_1, \psi_1 ) -  J(q_2, \psi_2 )    &  = & Pr \D \bar q  + Pr Ra (\bar T_y - \bar S_y)  \; , \;(y, z) \in D, \label{nneqn2}\\
\bar T_t + J(T_1, \psi_1 ) - J(T_2, \psi_2 )      &  = & \D \bar T \; , \;\; (y, z) \in D,\label{nneqn3}\\
\bar S_t + J(S_1, \psi_1 ) - J(S_2, \psi_2 )       &  = & \D \bar S \; , \;\; (y, z) \in D,\label{nneqn4}
\end{eqnarray}
The corresponding boundary conditions are as follows.

On the whole boundary:
\begin{equation}
\bar\psi=0, \;\bar\Delta \psi = \bar q = 0.
\end{equation}

\begin{equation}
\bar\theta_y (0, t) = \bar\theta_y (1, t)  = 0.
\end{equation}

At top $z = 1$:
\begin{equation}
\bar T_z + \bar T|_{z=1} = \bar\theta;\;    \bar S_z = 0.
 \end{equation}

At bottom $z = 0$:
\begin{equation}
\bar T_z = \bar S_z=0  \;.
\end{equation}

At the lateral boundary $ y= 0, 1$:
\begin{equation}
\bar T_y  = \bar S_y =0.
\end{equation}

Similar to the discussion in \S 2 above, we have
$$\frac12\frac{d}{dt}\|\bar\theta\|^2 = - \|\bar\theta_y\|^2 - \|\bar\theta\|^2 $$

\be
 - \int_0^1r(y)|\bar\theta|^2 dy + \int_0^1 \gamma(y)T(y, 1)\bar\theta dy + \int_0^1\bar f\bar\theta dy,
 \ee
\be
\frac12\frac{d}{dt}\|\bar q\|^2 + \int_D(J(q_1, \psi_1 ) -  J(q_2, \psi_2 ))\bar q = - Pr\|\nabla\bar q\|^2 +
PrRa\int_{\Omega}(\bar T_y - \bar S_y)\bar q,
 \ee
 \be \frac12\frac{d}{dt}\|\bar T\|^2 + \int_D(J(T_1, \psi_1 ) -  J(T_2, \psi_2 ))\bar T = - \|\nabla\bar  T\|^2 + \int_0^1[\bar\theta - \bar T(y, 1)]\bar T(y, 1) dy ,
 \ee

 \be
\frac12\frac{d}{dt}\|\bar S\|^2 + \int_D(J(S_1, \psi_1 ) -  J(S_2, \psi_2 ))\bar S = - \|\nabla\bar S\|^2.
 \ee
 In order to estimate the terms
 $\int_D(J(q_1, \psi_1 ) -  J(q_2, \psi_2 ))\bar q, \;\; \int_D(J(T_1, \psi_1 ) -  J(T_2, \psi_2 ))\bar T$
 and
  $\int_D(J(S_1, \psi_1 ) -  J(S_2, \psi_2 ))\bar S$,
  we need the following lemma:

  \begin{lemma}\label{lemma}
  The nonlinear Jacobian operaror $J(u, v)$ has the following property
  $$
 \|J(u_1, u_2 ) -  J(v_1, v_2 ))\| \le $$
 \be
 (\|\nabla u_1\| + \|\nabla u_2\| + \|\nabla v_1\|
 + \|\nabla v_2\|)(\|\nabla(u_1 - v_1)\| + \|\nabla(u_2 - v_2)\|),
 \ee
  for every $u_i, v_i \in H^1(D)$ (i = 1, 2).
  \end{lemma}

  The proof of this lemma is in \cite{gaoduan1}.

  By {\bf Lemma 3.1}, we have
  $$\int_D(J(q_1, \psi_1 ) -  J(q_2, \psi_2 ))\bar q \le (1 + \lambda_1)^2(\|\nabla q_1\| + \|\nabla q_2\|)\|\nabla\bar q\|\|\bar q\|,$$
  $$ \int_D(J(T_1, \psi_1 ) -  J(T_2, \psi_2 ))\bar T
  $$
  $$\le (\|\nabla T_1\| + \|\nabla T_2\| + \lambda_1(\|\nabla q_1\| + \|\nabla q_2\|))
  (\|\nabla\bar T\| + \lambda_1\|\nabla\bar q\|)\|\bar T\|,
  $$
  and
  $$\int_D(J(S_1, \psi_1 ) -  J(S_2, \psi_2 ))\bar S
  $$
  $$\le (\|\nabla S_1\| + \|\nabla S_2\| + \lambda_1(\|\nabla q_1\| + \|\nabla q_2\|))
  (\|\nabla\bar S\|+ \lambda_1\|\nabla\bar q\|)\|\bar S\|.
  $$
 Note that
 $$
 PrRa\int_{\Omega}(\bar T_y - \bar S_y)\bar q \le \frac{Pr}{2\lambda_1}\|\bar q\|^2 +
 {PrRa^2}{\lambda_1}(\|\bar T_y\|^2 + \|\bar S_y\|^2)$$
 $$
  \le \frac{Pr}{2}\|\nabla \bar q\|^2 +
 {PrRa^2}{\lambda_1}(\|\nabla\bar T\|^2 + \|\nabla \bar S\|^2).
 $$
Then by a similar argument as in \S 2 and using the Cauchy-Schwarz inequality,
we conclude that
 $$
 \frac12\frac{d}{dt}(\|\bar\theta\|^2 + \|\bar q\|^2 + 2{PrRa^2}{\lambda_1}\|\bar T\|^2 + 2{PrRa^2}{\lambda_1}\|\bar S\|^2) $$
 \be \le C_4(\|\bar\theta\|^2 + \|\bar q\|^2 + 2{PrRa^2}{\lambda_1}\|\bar T\|^2 + 2{PrRa^2}{\lambda_1}\|\bar S\|^2) + \|\bar f\|^2,
 \ee
 where $C_4$ depends on
 $Pr, Ra$, $\lambda_1$, $\|\gamma\|_{L^{\infty}}$
 as well as the $H^1-$norm of $q, T$ and $S$.
 By the Gronwall inequality, we further have
 \be
 \|\bar\theta\|^2 + \|\bar q\|^2 + \|\bar T\|^2 + \|\bar S\|^2
     \le C(b)\|\bar f\|^2,
 \ee
 for any given $b \; (0 < b < \infty, 0 < t \le b)$ and some positive constant $C(b)$
 depending on $b$ and $C_4$.
 Furthermore, we can obtain the similar estimates for the gradient
 of $\{\theta, q, T, S\}$, we omit the proof here, as the similar derivation
  will  be done in \S 5. Thus the solution differences
   $\bar\theta, \bar q, \bar T, \bar S$ and $\bar f$ are bounded when
   the external fluctuation difference $\bar f$ is bounded.
   So we have the following stability theorem.

 \begin{theorem}  \label{stability}
 ({\bf Stability under the external fluctuation})
 The  coupled atmosphere-ocean system
 (\ref{neqn1})--(\ref{neqn4}) (that is (\ref{eqn1})--(\ref{eqn4}))
 is stable under the external fluctuation in the atmospheric
 energy balance model. Namely, the solution of the coupled system
 depend continuously on the external fluctuation $f$ in $H^1$.
 \end{theorem}

%%%%%%%%%%%%%%%%%
%%%%%%%%%%%%%%%
%%%%%%%%%%%%%

\section{Dissipativity and Atmospheric Temperature Feedback}
\setcounter{equation}{0}

The ocean and the atmosphere are constantly interacting through
the air-sea exchange process. It is expected that the thermohaline circulation
could provide feedback to the air temperature. This is a direct impact of
 the thermohaline circulation on the climate. It is desirable to predict
or estimate this feedback.

To this end, let us estimate the  air temperature  $\theta$ in the mean-square
 norm,
in terms of the freshwater flux $F(y)$, external fluctuation $f(y,t)$ in the
energy balance model, the earth's longwave radiative cooling coefficient
$a$,
and  the empirical functions $S_a(y)$ and $ S_o(y)$
representing the latitudinal dependence of the shortwave solar radiation,
as well as physical parameters $Pr$ and $Ra$.

We will also show that the system generated by
(\ref{eqn1})--(\ref{side}) is a dissipative system in the sense of
\cite{hale} or \cite{temam} under some conditions, that is all
solutions $\{\theta, q, T, S\}$ enter a bounded set (so-called
{\em absorbing set}) in $H^1(0, 1)\times H_0^1(D)\times
H^1(D)\times H^1(D)$ after a finite time. Since
$$\frac{d}{dt}\int_{\Omega} S dydz = \int_0^1 F(y) dy = \mbox{constant}.
$$
For simplicity, we assume that
\be\label{mean}
\int_0^1 F(y) dy = 0, \;\;\int_DSdydz = 0
\ee
and
\be\label{gamma}
0 < \gamma(y) \le 1\; \mbox{or}\; 0 \le \gamma(y) < 1.
\ee
First, we derive a uniform estimate for $\{\theta, q, T, S\}$
in $L^2(0, 1)\times L^2(D) \times L^2(D) \times L^2(D)$.

Using the standard energy estimate as given in \S 2, we get

$$\frac12\frac{d}{dt}\|\theta\|^2 = - \|\theta_y\|^2 -
\|\theta\|^2 $$

$$
- a\int_0^1\theta dy + \int_0^1S_a\theta dy -
\int_0^1r(y)[S_o + \theta - T(y,1)]\theta dy + \int_0^1f\theta dy
 $$
 $$
 \le - \|\theta_y\|^2 - ( 1 - \e)\|\theta\|^2 - \inf\limits_{y\in[0,1]}\gamma(y)\|\theta\|^2
 $$
 \be\label{abtheta}
  + \frac{1}{\e}[a^2 + \|S_o\|^2 + \|S_a\|^2 + \|f\|^2] +  \|r(y)\|_{L^{\infty}}\int_0^1|\theta||T(y,1)|dy .
 \ee

 $$ \frac12\frac{d}{dt}(\|q\|^2 + 2PrRa^2\lambda_1(\| T\|^2 + \| S\|^2)) \le$$
 $$-\frac{Pr}{2}\|\nabla q\|^2  - PrPa^2\lambda_1(\|\nabla T\|^2 + \|\nabla S\|^2)
 $$
$$
+ 2PrRa^2\lambda_1\int_0^1[S_o + \theta - T(y,1)]T(y, 1) dy + 2PrRa^2\lambda_1\int_0^1F(y)S(y,1) dy.
 $$
 $$
 \le -\frac{Pr}{2}\|\nabla q\|^2  - PrPa^2\lambda_1(\|\nabla T\|^2 + \|\nabla S\|^2) - 2PrRa^2\lambda_1(1-\e)\int_0^1|T(y,1)|^2dy
 $$
 \be\label{abqts}
 + 2PrRa^2\lambda_1\int_0^1|\theta||T(y, 1)|dy + \frac{PrRa^2\lambda_1}{\e}\|S_o\|^2 + \frac{Pr^2Ra^4\lambda_1^2}{\e_1}\|F\|^2 + \e_1\|S(y, 1)\|^2,
 \ee
here $\e>0$ and $\e_1>0$ will be chosen later. By the trace
inequality, we have
$$
\|S(y, 1)\|^2 \le C(\|\nabla S\|^2 + \|S\|^2) \le C(1 + \bar\lambda_1)\|\nabla S\|^2,
$$
where $\bar\lambda_1$  is the constant in
the following Poincar\'e inequality (note that $\int_D S dydz = 0$)
\be\label{poincare}
\|S\|^2 \le \bar\lambda_1\|\nabla S\|^2.
\ee
Choosing $\e_1 = \frac{PrRa^2\lambda_1}{2(1+\bar\lambda_1)C}$, then (\ref{abqts}) can be written as
$$ \frac12\frac{d}{dt}(\|q\|^2 + 2PrRa^2\lambda_1(\| T\|^2 + \| S\|^2)) \le$$
 $$-\frac{Pr}{2}\|\nabla q\|^2  - PrPa^2\lambda_1(\|\nabla T\|^2 + \frac12\|\nabla S\|^2) - 2PrRa^2\lambda_1(1-\e)\int_0^1|T(y,1)|^2dy
 $$
  \be\label{abqtsineq}
+ 2PrRa^2\lambda_1\int_0^1|\theta||T(y, 1)|dy + \frac{PrRa^2\lambda_1}{\e}\|S_o\|^2 + 2PrRa^2\lambda_1(1 + \bar\lambda_1)C\|F\|^2.
 \ee
 Then, multiplying (\ref{abtheta}) by $2PrRa^2\lambda_1$ and
 adding to (\ref{abqtsineq}), we get

 $$ \frac12\frac{d}{dt}(\|q\|^2 + 2PrRa^2\lambda_1(\|\theta\|^2 + \| T\|^2 + \| S\|^2)) \le$$
 $$- 2PrRa^2\lambda_1\|\theta_y\|^2 -\frac{Pr}{2}\|\nabla q\|^2  - PrPa^2\lambda_1(\|\nabla T\|^2 + \frac12\|\nabla S\|^2)
 $$
 $$- 2PrRa^2\lambda_1((1-\e)\|\theta\|^2 - \inf\limits_{y\in[0,1]}\gamma(y)\|\theta\|^2 - (1 + \|\gamma\|_{L^{\infty}})
 \int_0^1|T(y,1)||\theta| dy + (1-\e)\int_0^1|T(y,1)|^2dy)$$
   \be\label{abtotal}
 + \frac{PrRa^2\lambda_1}{2\e}[a^2 + \frac54\|S_o\|^2 + + \|S_a\|^2 + \|f\|^2] + 2PrRa^2\lambda_1(1 + \bar\lambda_1)C\|F\|^2.
 \ee
 By (\ref{gamma}), when $0 \le \gamma < 1$,
 we could choose $\e$ such that
 (since $\inf\limits_{y\in[0,1]} \gamma(y) = 0$ now)

 $$4(1-\e)^2 > (1 + \|\gamma\|_{L^{\infty}})^2, \;\mbox{i.e.}\;
 \e < \frac{1 - \|\gamma\|_{L^{\infty}}}{1 + \|\gamma\|_{L^{\infty}}}
 := \alpha_0.$$
 For example, we choose $\e = \frac{\alpha_0}{2}$, then
 $$- 2PrRa^2\lambda_1((1-\e)\|\theta\| - (1 + \|\gamma\|_{L^{\infty}})
 \int_0^1|T(y,1)||\theta| dy + (1-\e)\int_0^1|T(y,1)|^2dy)$$
 $$ < - \frac{PrRa^2\lambda_1\alpha_0}{2}(\|\theta\|^2 + \|T(y, 1)\|^2).$$

 If $0 < \gamma(y) \le 1$ as in (\ref{gamma}),
 we denote $\inf\limits_{y\in[0,1]}\gamma(y) = \beta_0$.
 Then we take $\e = \frac{\beta_0}{6}$ to obtain
 $$- 2PrRa^2\lambda_1((1-\e + \beta_0)\|\theta\| - (1 + \|\gamma\|_{L^{\infty}})
 \int_0^1|T(y,1)||\theta| dy + (1-\e)\int_0^1|T(y,1)|^2dy)$$
 $$ < - \frac{PrRa^2\lambda_1\beta_0}{6}(\|\theta\|^2 + \|T(y, 1)\|^2).$$

 So, in the case of $0 \le \gamma(y) < 1$, (\ref{abtotal}) can be written as

 $$ \frac12\frac{d}{dt}(\|q\|^2 + 2PrRa^2\lambda_1(\|\theta\|^2 + \| T\|^2 + \| S\|^2)) \le$$
 $$- 2PrRa^2\lambda_1\|\theta_y\|^2 -\frac{Pr}{2}\|\nabla q\|^2  - PrRa^2\lambda_1(\|\nabla T\|^2 + \frac12\|\nabla S\|^2)
 $$
 $$- \frac{PrRa^2\lambda_1\alpha_0}{2}(\|\theta\|^2 + \|T(y, 1)\|^2)$$
   \be\label{abtotal1}
 + \frac{PrRa^2\lambda_1}{\alpha_0}[a^2 + \frac54\|S_o\|^2 + \|S_a\|^2 + \|f\|^2] + 2PrRa^2\lambda_1(1 + \bar\lambda_1)C\|F\|^2.
 \ee
 For $0 < \gamma(y) < 1$, we will have similar estimate. Since
$$ T^2(y,z) - T^2(y,1) = 2\int_1^z TT_z dz,$$
  we further have

\be\label{tcoest} \|T\|^2 \le 2\int_0^1|T(y,1)|^2dy + 4\|\nabla T
\|^2. \ee

Using the Poincar\'e inequlity again for $q$ and
letting $\alpha_1 = \min\{\frac{Pr\lambda_1}{4}, \frac{\alpha_0}{8}, \frac{1}{8},\frac{\bar\lambda_1}{8}\}$ and $\beta =
\min\{\frac{Pr\lambda_1}{4}, $
$\frac{PrRa^2\lambda_1\alpha_0}{4}, \frac{PrRa^2\lambda_1}{4}\}$,
the estimate (4.8) becomes

 $$ \frac12\frac{d}{dt}(\|q\|^2 + 2PrRa^2\lambda_1(\|\theta\|^2 + \| T\|^2 + \| S\|^2)) \le$$
 $$- \alpha_1(\|q\|^2 + 2PrRa^2\lambda_1(\|\theta\|^2 + \|T\|^2 + \|S\|^2) $$
 $$- \beta_1(\|\theta_y\|^2 + \|\nabla q\|^2  + \|\nabla T\|^2 + \|\nabla S\|^2 + \|T(y, 1)\|^2)
 $$
    \be\label{total2}
 + \frac{PrRa^2\lambda_1}{\alpha_0}[a^2 + \frac54\|S_o\|^2 + \|f\|^2] + 2PrRa^2\lambda_1(1 + \bar\lambda_1)C\|F\|^2.
 \ee
Using the Gronwall inequality, we finally obtain the mean-square norm
estimate for the solution of the coupled atmosphere-ocean model:
%%%%%%%%%%%%%%%%%
%%%%%%%%%%%%%%%%%
%%%%%%%%%%%%%%%%%
$$\|q\|^2 + 2PrRa^2\lambda_1(\|\theta\|^2 + \| T\|^2 + \| S\|^2) \le $$
$$
e^{-\alpha_1t}(\|q_0\|^2 + 2PrRa^2\lambda_1(\|\theta_0\|^2 + \| T_0\|^2 + \| S_0\|^2)) $$
\be\label{abl2}
 + \frac{PrRa^2\lambda_1}{\alpha_1}[\frac{1}{\alpha_0}(a^2 + \frac54\|S_o\|^2 + \|S_a\|^2 + \|f\|^2)
+ 2(1 + \bar\lambda_1)C\|F\|^2]
\ee

\bigskip

In particular, we get the mean-square norm estimate for the atmospheric
temperature feedback
$$
\|\theta\|^2
\le
e^{-\alpha_1t}(\frac1{2PrRa^2\lambda_1} \|q_0\|^2
+ \|\theta_0\|^2 + \| T_0\|^2 + \| S_0\|^2 )
$$
\be\label{feedback}
 + \frac{1}{2\alpha_1}[\frac{1}{\alpha_0}(a^2 + \frac54\|S_o\|^2 + \|S_a\|^2 + \|f\|^2)
+ 2(1 + \bar\lambda_1)C\|F\|^2].
\ee
This atmospheric temperature feedback estimate is in terms of
physical quantities such as the freshwater flux $F(y)$, external fluctuation $f(y,t)$ in the
energy balance model, the earth's longwave radiative cooling coefficient
$a$,
and  the empirical functions $S_a(y)$ and $ S_o(y)$
representing the latitudinal dependence of the shortwave solar radiation,
as well as physical parameters $Pr$ and $Ra$. Here $\lambda_1$ and
$ \bar\lambda_1 $ are the constants in the Poincar\'e  inequality on the
 domain $D$
in the cases of zero Dirichlet boundary condition and zero mean value,
respectively. Moreover,  $C$ is a constant depending only on the domain $D$,
$\alpha_0=\frac{1 - \|\gamma\|_{L^{\infty}}}{1 + \|\gamma\|_{L^{\infty}}}$,
and
$\alpha_1 = \min\{\frac{Pr\lambda_1}{4}, \frac{\alpha_0}{8},
\frac{1}{8},\frac{\bar\lambda_1}{8}\}$.

\bigskip

We can furthermore derive solution estimate in $H^1$ norm. To do so, we  first
 get from (\ref{total2})
$$
\int_t^{t+1}(\|\theta_y\|^2 + \|\nabla q\|^2
+ \|\nabla T\|^2 + \|\nabla S\|^2 + \|T(y, 1)\|^2) \le$$
$$\frac{1}{\beta_1}e^{-\alpha_0t}(\|q_0\|^2
+ 2PrRa^2\lambda_1(\|\theta_0\|^2 + \| T_0\|^2 + \| S_0\|^2)) $$
\be\label{abl20t}
 + \frac{2PrRa^2\lambda_1}{\alpha_1\beta_1}(\frac{1}{\alpha_0}(a^2
 + \frac54\|S_o\|^2 + \|S_a\|^2 + \|f\|^2)
+ 2(1 + \bar\lambda_1)C\|F\|^2).
\ee

So, let $\|q_0\|^2 + 2PrRa^2\lambda_1(\|\theta_0\|^2
+ \| T_0\|^2 + \| S_0\|^2)$ be bounded by
some (big) upper bound $R^2$ and denote $M^2 =
\frac{PrRa^2\lambda_1}{\alpha_1}(\frac{1}{\alpha_0}(a^2 + \frac54\|S_o\|^2
+ \|S_a\|^2 + \|f\|^2)
+ 2(1 + \bar\lambda_1)C\|F\|^2)$. Then there is a time
$t^* \ge \frac{2}{\alpha_1}\ln{\frac{R}{M}}$ such that
\be \label{abballl2}
\|q\|^2 + 2PrRa^2\lambda_1(\|\theta\|^2 + \| T\|^2
+ \| S\|^2) \le 2M^2,\;t \ge t^*
\ee
and
\be\label{abballl20t}
\int_t^{t+1}(\|\theta_y\|^2 + \|\nabla q\|^2
+ \|\nabla T\|^2 + \|\nabla S\|^2
+ \|T(y, 1)\|^2) \le \frac{3}{\beta_1}M^2, \; t \ge t^*.
\ee
Next, we derive a uniform estimate of gradient of $\{\theta, q, T, S\}$
in $L^2(0, 1)\times L^2(D) \times L^2(D) \times L^2(D)$.
In order to avoid the difficulty caused by the
non-homogeneous boundary conditions,
we use equations (\ref{neqn1})--(\ref{side1})
 instead of (\ref{eqn1})--(\ref{side}).

Multiplying (\ref{neqn1})--(\ref{neqn4})
by $-\theta_{yy}, -\D q, -\D T$ and $-\D S$
respectively, integrating over $(0, 1)$ and $D$,
noting that $S^*\in H^2(D)$ is
known and $T^*$ is only dependent on
$\theta$, we get
$$\frac12\frac{d}{dt}\|\theta_y \|^2 = - \|\theta_{yy}\|^2 - \|\theta_y\|^2 $$
\be \label{dtheta}
+  a\int_0^1\theta_{yy} dy -  \int_0^1S_a\theta_{yy} dy +
\int_0^1r(y)[S_o + \theta - T(y,1)]\theta_{yy} dy - \int_0^1f\theta_{yy} dy,
 \ee
 \be\label{dq}
 \frac12\frac{d}{dt}\|\nabla q\|^2  = - Pr\|\D q\|^2 +  \int_DJ(q, \psi)\D q - \int_DPrRa(T_y - S_y + T_y^* - S_y^*)\D q,
 \ee
 $$
 \frac12\frac{d}{dt}(\|\nabla T\|^2 + \|T(y, 1)\|^2) = - \|\D T\|^2  + \int_DJ(T + T^*, \psi)\D T
 $$
 $$ - \int_D[((1 - e^{1-z})(1 - \gamma(y)) - e^{1-z})\eta_{\e}(z) + (1 - e^{1-z})\eta''_{\e}(z) + 2e^{1-z}\eta'_{\e}(z) - e^{1-z}\eta_{\e}]\theta\D T$$
\be\label{dt}
 + \int_D (1 - e^{1-z})\eta_{\e}(z)\gamma(y)T(y, 1)\D T - \int_Dg\D T,
 \ee
 $$
 \frac12\frac{d}{dt}\|\nabla S\|^2 = - \|\D S\|^2  + \int_DJ(S + S^*, \psi)\D S
 $$
 \be\label{ds}
 - \int_DF''(y)z\eta_{\e}(z)\D S - \int_DF(y)(2\eta_{\e}^{\prime} + z\eta_{\e}{''}(z)) \D S.
  \ee
  Note that
  $$
  a\int_0^1\theta_{yy} dy -  \int_0^1S_a\theta_{yy} dy +
\int_0^1r(y)[S_o + \theta - T(y,1)]\theta_{yy} dy - \int_0^1f\theta_{yy} dy
$$
\be\label{dthetaineq}
\le \frac{3\e}{2}\|\theta_{yy}\|^2 + \frac{1}{4\e}[a^2 + \|S_a\|^2 + \|S_o\|^2 + \|f\|^2 + \|\theta\|^2 + \|T(y, 1)\|^2],
\ee
$$
- \int_DPrRa(T_y - S_y + T_y^* - S_y^*)\D q \le
$$
\be\label{dts}
\frac{Pr}{2}\|\D q\|^2 + \frac{5PrRa^2}{2}(\|\nabla T\|^2 + \|\nabla S\|^2 + \|\theta_y\|^2)
+ \frac{5PrRa^2}{2}(\|S_o'\|^2 + \|F'\|^2),
\ee
 $$ - \int_D[((1 - e^{1-z})(1 - \gamma(y)) - e^{1-z})\eta_{\e}(z) + (1 - e^{1-z})\eta''_{\e}(z) + 2e^{1-z}\eta'_{\e}(z) - e^{1-z}\eta_{\e}]\theta\D T$$
$$ + \int_D (1 - e^{1-z})\eta_{\e}(z)\gamma(y)T(y, 1)\D T - \int_Dg\D T
 $$
 \be\label{dtn}
 \le \frac12\|\D T\|^2 + C(a^2 + \|S_a\|^2 + \|S_o\|_{H^2}^2 + \|f\|^2 + \|\theta\|^2 + \|T(y, 1)\|^2).
 \ee
 \be\label{dsn}
 - \int_DF''(y)z\eta_{\e}(z)\D S - \int_DF(y)(2\eta_{\e}^{\prime} + z\eta_{\e}{''}(z)) \D S \le
  \frac12\|\D S\|^2 + C\|F\|_{H^2}^2.
   \ee
   About the estimates of $\int_DJ(q, \psi)\D q, \int_DJ(T + T^*, \psi)\D T$ and $\int_DJ(S + S^*, \psi)\D S$, similar to \cite{gaoduan}, we have
  \be\label{dqj}
   \int_DJ(q, \psi)\D q \le C\|\D\psi\|\|\nabla q\|\|\D q\| = C\|q\|\|\nabla q\|\|\D q\|,
   \ee
   \be\label{dtj}
   \int_DJ(T + T^*, \psi)\D T \le C\|q\|(\|\D T\| + \|\nabla T\| + \|T\|)\|\nabla T\| + \|q\|(\|\theta_y\| + \|S_o'\|)\|\D T\|,
    \ee
   \be\label{dsj}
   \int_DJ(S + S^*, \psi)\D S \le C\|q\|(\|\D S\| + \|\nabla S\| + \|S\|)\|\nabla S\| + \|q\|\|F'\|\|\D S\|.
    \ee
    Using the Cauchy-Schwarz inequality, (\ref{dthetaineq})--(\ref{dsj})
    and (\ref{abballl2}), when $t \ge t^*$, we get
 $$ \frac12\frac{d}{dt}(\|\theta_y\|^2 + \|\nabla q\|^2 + 3PrRa^2\lambda_1(\|\nabla T\|^2 + \|T(y, 1)\|^2 + \|\nabla S\|^2)) $$
 $$\le C(\|\theta_y\|^2 + \|\nabla q\|^2 + 3PrRa^2\lambda_1(\|\nabla T\|^2 + \|T(y, 1)\|^2 + \|\nabla S\|^2)
 $$
 \be\label{dtotal}
 + C(a^2 + \|S_a\|^2 + \|S_o\|_{H^2}^2 + \|f\|^2 + \|F\|_{H^2}^2).
  \ee
  By (\ref{abballl20t}) and a uniform Gronwall lemma (\cite{temam}), we obtain
 $$
  \|\theta_y\|^2 + \|\nabla q\|^2 + 3PrRa^2\lambda_1(\|\nabla T\|^2 + \|T(y, 1)\|^2 + \|\nabla S\|^2)
  $$
  \be
  \label{dtotalineq}
   \le C(a^2 + \|S_a\|^2 + \|S_o\|_{H^2}^2 + \|f\|^2 + \|F\|_{H^2}^2).
  \ee
 By (\ref{abballl2}) and (\ref{dtotalineq}), we know there exists an
  absorbing sets ${\cal B}$ in
  $H^1(0, 1)\times H^1_0(D)\times H^1(D)\times H^1(D)$ for the solution of
 (\ref{eqn1})--(\ref{side}):
 \be\label{ball}
   {\cal B} = {\large\{}\{\theta, q, T, S\}: \;\;
  \|\theta\|_1^2  + \|q\|_1^2 + \|T\|_1^2 + \|S\|_1^2 \le
 C(a^2 + \|S_a\|^2 + \|S_o\|_{H^2}^2 + \|f\|^2 + \|F\|_{H^2}^2){\large\}},
 \ee
 that is, for every bounded set in
 $H^1(0, 1)\times H^1(D)\times H^1(D)\times H^1(D)$, when
 $t \ge t^* + 1$, the solution of (\ref{eqn1})--(\ref{side})
 will enter into the ${\cal B}$.

We summarize our results in section in the following theorem.

\begin{theorem}\label{feedbacktheorem}
({\bf Atmospheric temperature feedback and dissipativity})
Assume that the freshwater flux $F(y)$ has zero mean as in  (\ref{mean})
and the ocean basin's latitudinal fraction function $\gamma(y)$
is bounded as in (\ref{gamma}).
Then the coupled atmosphere-ocean system
 (\ref{eqn1})--(\ref{eqn4}) has an absorbing set  in
 $H^1(0, 1)\times H^1_0(D)\times H^1(D)\times H^1(D)$
as given in (\ref{ball}).
 More importantly, the atmospheric temperature feedback $\theta(y,t)$
 is bounded in mean-square norm in terms of
physical quantities such as the freshwater flux $F(y)$,
external fluctuation $f(y,t)$ in the
energy balance model, the earth's longwave radiative cooling coefficient
$a$,
and  the empirical functions $S_a(y)$ and $ S_o(y)$
representing the latitudinal dependence of the shortwave solar radiation,
as well as the  Prandtl  number $Pr$ and  the Rayleigh number $Ra$
as in (\ref{feedback}).
\end{theorem}

\begin{remark}
Due to $\frac{d}{dt}\int_D \bar S = 0$, and $\int_D\bar S_0 = 0$,
 we obtain $\int_D\bar S = 0$. So, as seen in
the discussion of this section, we know that  the
stability (proved in \S 3) under the exernal fluctuation is uniform  in time $t$ when $0 < \gamma(y) \le 1$ or $0 \le \gamma(y) < 1$.
\end{remark}

%%%%%%%%%%%%%%%%%%%%%%%%%%%%%%
%%%%%%%%%%%%%%%%%%%%%%%%%%%%%
%%%%%%%%%%%%%%%%%%%
\section{Strong Contraction and Almost Periodic Atmosphere-Ocean Dynamics}
\setcounter{equation}{0}

In this section, we study the coupled atmosphere-ocean dynamical response
to almost periodic (in particular, periodic and quasi-periodic)
external fluctuation $f(y,t)$ in the atmospheric energy
balance model (\ref{eqn1}). A central question is: Does the
coupled atmosphere-ocean system respond almost
periodically  to almost  periodic external fluctuation $f(y,t)$?

To answer this question, we need to understand
the strong contraction property of the
coupled atmosphere-ocean  system in the absorbing set
${\cal B}$ defined in (\ref{ball}). Let $\{\theta^i, q^i, T^i, S^i\}$ be
 two trajectories corresponding to initial values
$\{\theta_0^i, q_0^i, T_0^i, S_0^i\}\in {\cal B}$ for $ i = 1, 2$.
Note that these trajectories remain inside ${\cal B}$
as $B$ is a forward invariant set. Their difference
$$\de\theta = \theta^1 - \theta^2,\;\de q = q^1 - q^2,\;\de T = T^1 - T^2,\;\de S = S^1 - S^2$$
satisfy the following equations:
\begin{eqnarray}
\de\theta_t   =   \de\theta_{yy} - \de\theta &-&\gamma(y)[\de\theta - \de T], \;\;\;\;\;\;\;\; 0\le z \le 1,\label{cneqn1}\\
\de q_t + J(q_1, \psi_1 ) -  J(q_2, \psi_2 )    &  = & Pr \D \de q  + PrRa (\de T_y - \de S_y)  \; , \;(y, z) \in D, \label{cneqn2}\\
\de T_t + J(T_1, \psi_1 ) - J(T_2, \psi_2 )      &  = & \D \de T \; , \;(y, z) \in D,\label{cneqn3}\\
\de S_t + J(S_1, \psi_1 ) - J(S_2, \psi_2 )       &  = & \D \de S \; , \;(y, z) \in D,\label{cneqn4}
\end{eqnarray}
The corresponding boundary conditions are:
\begin{equation}
\de \theta_y (0, t) = \de\theta_y (1, t)  = 0.
\end{equation}

On the whole boundary:
\begin{equation}
\de\psi=0, \; \de q = 0.
\end{equation}

At top $z = 1$:
\begin{equation}
\de T_z + \de T|_{z=1} = \de\theta;\;    \de S_z = 0.
 \end{equation}

At bottom $z = 0$:
\begin{equation}
\de T_z = \de S_z=0.
\end{equation}

At the lateral boundary $ y= 0, 1$:
\begin{equation}
\de T_y  = \de S_y =0.
\end{equation}
The initial   conditions are:
\be
\de\theta_0 = \theta^1_0 - \theta^2_0,\; \de q_0 = q_0^1 - q_0^2,\; \de T_0 = T_0^1 - T_0^2,\; \de S_0 = S_0^1 - S_0^2,
\ee
where $\de\theta_0 = \de\theta(y, 0), \de q_0 = \de q(y, z, 0), \de T_0 = \de T(y, z, 0)$ and $\de S_0 = \de S(y, z, 0)$.

 Using energy estimates as in \S 3, we have
\be\label{ctheta}
\frac12\frac{d}{dt}\|\de\theta\|^2 = - \|\de\theta_{y}\|^2 - \|\de\theta\|^2 - \int_D\gamma(y)|\de\theta|^2 dy- \int_D \gamma(y)\de\theta\de T(y,1) dy,
\ee
 \be\label{cq}
 \frac12\frac{d}{dt}\|\de q\|^2  = - Pr\|\nabla \de q\|^2 -  \int_D(J(q^1, \psi^1) - J(q^2, \psi^2))\de q - \int_DPrRa(\de T_y - \de S_y)\de q,
 \ee
 $$
 \frac12\frac{d}{dt}\|\de T\|^2 = - \|\nabla\de T\|^2
 $$
 \be\label{ct}
  + \int_0^1[\de\theta - \de T(y, 1)]\de T(y, 1) dy - \int_D(J(T^1, \psi^1) - J(T^2, \psi^2))\de T,
 \ee
 \be\label{cs}
 \frac12\frac{d}{dt}\|\de S\|^2 = - \|\nabla\de S\|^2  - \int_D(J(S^1, \psi^1) - J(S^2, \psi^2))\de S.
 \ee

 Using {\bf Lemma \ref{lemma}}, we imply
 $$
 -  \int_D(J(q^1, \psi^1) - J(q^2, \psi^2))\de q \le \|J(q^1, \psi^1) - J(q^2, \psi^2)\|\|\de q\| \le
 $$
 $$(\|\nabla q^1\| + \|\nabla q^2\| + \|\nabla\psi^1\| + \|\nabla\psi^2\|)(\|\nabla\de q\| + \|\nabla\de\psi\|)\|\de q\| $$
 $$\le (1 + {\lambda_1})^2\sqrt{\lambda_1}(\|q^1\| + \|q^2\|)\|\nabla\de q\|^2,
 $$
 $$
 - \int_D(J(T^1, \psi^1) - J(T^2, \psi^2))\de T \le  \|J(T^1, \psi^1) - J(T^2, \psi^2)\|\|\de T\|$$
 $$\le (\|\nabla T^1\| + \|\nabla T^2\| + \|\nabla\psi^1\| + \|\nabla\psi^2\|)(\|\nabla\de T\| + \|\nabla\de\psi\|)\|\de T\|
 $$
 $$\le (\|\nabla T^1\| + \|\nabla T^2\| + \lambda_1(\|\nabla q^1\| + \|\nabla q^2\|))((\|\nabla\de T\| + \lambda_1\|\nabla\de q\|)\|\de T\|,
 $$
 $$
 - \int_D(J(S^1, \psi^1) - J(S^2, \psi^2))\de S \le \|J(S^1, \psi^1) - J(S^2, \psi^2)\|\|\de S\|$$
 $$\le (\|\nabla S^1\| + \|\nabla S^2\| + \|\nabla\psi^1\| + \|\nabla\psi^2\|)(\|\nabla\de S\| + \|\nabla\de\psi\|)\|\de S\|
 $$
 $$\le (\|\nabla S^1\| + \|\nabla S^2\| + \lambda_1(\|\nabla q^1\| + \|\nabla q^2\|))((\|\nabla\de S\| + \lambda_1\|\nabla\de q\|)\|\de S\|.
 $$
 Other terms in (\ref{ctheta})--(\ref{cs}) can be estimated as in
 the proof of the existence of the absorbing set in the last section.
 So we have
 $$
\frac12\frac{d}{dt}(\|\de q\|^2 + 2PrRa^2\lambda_1(\|\de\theta\|^2 + \|\de T\|^2 + \|\de S\|^2)) \le
$$
$$
- \alpha_1(\|\de q\|^2 + 2PrRa^2\lambda_1(\|\de\theta\|^2 + \|\de T\|^2 + \|\de S\|^2)) $$
 $$- \beta_1(\|\de\theta_y\|^2 + \|\nabla\de q\|^2  + \|\nabla\de T\|^2 + \frac12\|\nabla\de S\|^2 + \|\de T(y, 1)\|^2)
$$
$$
+ (1 + {\lambda_1})^2\sqrt{\lambda_1}(\|q^1\| + \|q^2\|)\|\nabla\de q\|^2 + $$
$$
2PrRa^2\lambda_1(\|\nabla T^1\| + \|\nabla T^2\| + \lambda_1(\|\nabla q^1\| + \|\nabla q^2\|))((\|\nabla\de T\| + \lambda_1\|\nabla\de q\|)\|\de T\|$$
$$
 + 2PrRa^2\lambda_1(\|\nabla S^1\| + \|\nabla S^2\| + \lambda_1(\|\nabla q^1\| + \|\nabla q^2\|))((\|\nabla\de S\| + \lambda_1\|\nabla\de q\|)\|\de S\|.$$

Now we assume that $C(a^2 + \|S_a\|^2 + \|S_o\|_{H^2}^2
+ \|f\|^2 + \|F\|_{H^2}^2)$
is small enough. This is a condition imposed on the physical quantities
such as  the freshwater flux $F(y)$,
external fluctuation $f(y,t)$ in the
energy balance model, the earth's longwave radiative cooling coefficient
$a$,
and  the empirical functions $S_a(y)$ and $ S_o(y)$
representing the latitudinal dependence of the shortwave solar radiation.

 By the Cauchy-Schwarz inequality, (\ref{poincare}), (\ref{tcoest})
 and
  (\ref{ball}),  we get
$$
\frac12\frac{d}{dt}(\|\de q\|^2 + 2PrRa^2\lambda_1(\|\de\theta\|^2 + \|\de T\|^2 + \|\de S\|^2)) \le
$$
$$
- \alpha_1(\|\de q\|^2 + 2PrRa^2\lambda_1(\|\de\theta\|^2 + \|\de T\|^2 + \|\de S\|^2)) $$
 \be\label{ctotal1}
 - \frac{\beta_1}{2}(\|\de\theta_y\|^2 + \|\nabla\de q\|^2  + \|\nabla\de T\|^2 + \frac12\|\nabla\de S\|^2 + \|\de T(y, 1)\|^2).
\ee
By the Gronwall's inequality, we obtain the strong contraction
in $L^2(0, 1)\times L^2(D) \times L^2(D) \times L^2(D)$ for
the solution of the coupled atmosphere-ocean system
(\ref{eqn1})--(\ref{side}), that is
$$
\|\de q\|^2 + 2PrRa^2\lambda_1(\|\de\theta\|^2 + \|\de T\|^2 + \|\de S\|^2) \le
$$
\be\label{ctotal2}
 e^{-\alpha_1 t}(\|\de q_0\|^2 + 2PrRa^2\lambda_1(\|\de\theta_0\|^2 + \|\de T_0\|^2 + \|\de S_0\|^2)).
\ee
Next, we can show the strong contraction of gradient
 in $L^2(0, 1)\times L^2(D) \times L^2(D) \times L^2(D)$.
 We also need to estimate it using (\ref{neqn1})--(\ref{side1}).
 Noticing that $S^*$ is independent of $\{\theta, q, T, S\}$ and $T^*$
 is only dependent on $\theta$,    we get
 \be
\de\theta_t   =   \de\theta_{yy} - \de\theta  - \gamma(y)[\de\theta - \de T(y, 1)], \;  0\le z \le 1,\label{cnneqn1}
\ee
$$
\de q_t   =  Pr \D\de q  $$
\be
- J(q^1, \psi^1) +  J(q^2, \psi^2)+ PrRa (\de T_y - \de S_y) + \de\theta(1- e^{1-z})\eta_{\e}(z)),\;  (y, z) \in D, \label{cnneqn2}
\ee
$$
\de T_t  =  \D\de T - J(T^1, \psi^1) +  J(T^2, \psi^2) $$
$$- J((1-e^{1-z})\eta_{\e}(z)\theta^1, \psi^1) +  J((1-e^{1-z})\eta_{\e}(z)\theta^1, \psi^2)
- J((1-e^{1-z})\eta_{\e}(z)S_o(y), \de\psi)$$
$$
[((1 - e^{1-z})(1 - \gamma(y)) - e^{1-z})\eta_{\e}(z) + (1 - e^{1-z})\eta''_{\e}(z) + 2e^{1-z}\eta'_{\e}(z) - e^{1-z}\eta_{\e}]\de\theta$$
\be -
 (1 - e^{1-z})\eta_{\e}(z)\gamma(y)\de T(y, 1), \; (y, z) \in D,\label{cnneqn3}
 \ee
 \be
\de S_t = \D\de S - J(S^1, \psi^1) + J(S^2, \psi^2)  +   J(S_{\e}^*, \de\psi ), \; (y, z) \in D.\label{cnneqn4}
\ee

The corresponding boundary conditions are:
\begin{equation}
\de\theta_y (0, t) = \de\theta_y (1, t)  = 0.
\end{equation}

On the whole boundary:
\begin{equation}
\de\psi=0, \; \de\Delta \psi = \de q = 0.
\end{equation}

 At top $z=1$:
\begin{equation}
 \de T_z + \de T|_{z=1} = 0; \;\;   \de S_z = 0 .
  \end{equation}

At bottom $z=0$:
\begin{equation}
\de T_z =  \de S_z=0  \;.
\end{equation}

At  the lateral boundary $y=0, 1$:
\begin{equation}
\de T_y  = \de S_y =0.
\end{equation}

The initial   conditions are:
\be
\de\theta_0 = \theta^1_0 - \theta^2_0,\; \de q_0 = q_0^1 - q_0^2,\; \de T_0 = T_0^1 - T_0^2 + (1 - e^{1-z})\eta_{\e}(z)\de\theta_0,\; \de S_0 = S_0^1 - S_0^2,
\ee
where $\de\theta_0 = \de\theta(y, 0), \de q_0 = \de q(y, z, 0), \de T_0 = \de T(y, z, 0)$ and $\de S_0 = \de S(y, z, 0)$.

Multiplying (\ref{cnneqn1})--(\ref{cnneqn4}) by $-\de\theta_{yy}, -\D\de q, -\D\de T$ and $-\D\de S$ respectively, integrating over $(0, 1)$ and $D$,
 we get
\be\label{cdtheta}
\frac12\frac{d}{dt}\|\de\theta_y \|^2 = - \|\de\theta_{yy}\|^2 - \|\de\theta_y\|^2
+ \int_0^1r(y)[\de\theta - \de T(y,1)]\de\theta_{yy} dy,
 \ee
 $$
 \frac12\frac{d}{dt}\|\nabla\de q\|^2  = - Pr\|\D\de q\|^2 +  \int_D (J(q^1, \psi^1) -  J(q^2, \psi^2))\D\de q $$
 \be\label{cdq}
  - PrRa\int_D (\de T_y - \de S_y )\D\de q - \int_D(1- e^{1-z})\eta_{\e}(z))\de\theta\D\de q,
 \ee
 $$
 \frac12\frac{d}{dt}(\|\nabla\de T\|^2 + \|\de T(y, 1)\|^2) = - \|\D\de T\|^2  + \int_D(J(T^1, \psi^1) -  J(T^2, \psi^2))\D\de T
 $$
 $$+ \int_D(J((1-e^{1-z})\eta_{\e}(z)\theta^1, \psi^1) -  J((1-e^{1-z})\eta_{\e}(z)\theta^2, \psi^2))\D\de T $$
$$
- \int_D[((1 - e^{1-z})(1 - \gamma(y)) - e^{1-z})\eta_{\e}(z) + (1 - e^{1-z})\eta''_{\e}(z) + 2e^{1-z}\eta'_{\e}(z) - e^{1-z}\eta_{\e}]\de\theta\D\de T$$
$$
+ \int_D (1 - e^{1-z})\eta_{\e}(z)\gamma(y)\de T(y, 1)\D\de T
 $$
\be\label{cdt}
+ \int_DJ((1-e^{1-z})\eta_{\e}(z)S_o(y), \de\psi)\D\de T, \; (y, z) \in D,
 \ee

  $$
 \frac12\frac{d}{dt}\|\nabla\de S\|^2 = - \|\D\de S\|^2  $$
 \be\label{cds}
 + \int_D(J(S^1, \psi^1) - J(S^2, \psi^2))\D\de S +   \int_DJ(S_{\e}^*, \de\psi )\D\de S,\; (y, z) \in D.
 \ee

  Note that
  $$
  \int_0^1r(y)[\de\theta - \de T(y,1)]\de\theta_{yy} dy \le \frac12\|\de\theta_{yy}\|^2 + \|\gamma(y)\|_{L^{\infty}}\|^2(\|\de\theta\|^2 + \|\de T(y,1)\|^2),
$$
$$
- PrRa\int_D (\de T_y - \de S_y )\D\de q \le \frac{Pr}{4}\|\D\de q\|^2 + 2PrRa^2(\|\nabla\de T\|^2 + \|\nabla\de S\|^2),
$$
$$
- \int_D(1- e^{1-z})\eta_{\e}(z))\de\theta\D\de q \le \frac{Pr}{4}\|\D\de q\|^2 + \frac{1}{Pr}\|\de\theta\|^2,
$$
$$
- \int_D[((1 - e^{1-z})(1 - \gamma(y)) - e^{1-z})\eta_{\e}(z) + (1 - e^{1-z})\eta''_{\e}(z) + 2e^{1-z}\eta'_{\e}(z) - e^{1-z}\eta_{\e}]\de\theta\D\de T
$$
$$
\le \frac16\|\D\de T\|^2 + C\|\de\theta\|^2,
$$
$$
\int_D (1 - e^{1-z})\eta_{\e}(z)\gamma(y)\de T(y, 1)\D\de T \le \frac16\|\D\de T\|^2  + \frac{3}{2}\|\de T(y, 1)\|^2,
$$
$$
\int_DJ((1-e^{1-z})\eta_{\e}(z)S_o(y), \de\psi)\D\de T \le \frac16\|\D\de T\|^2 + C\lambda_1(\|S_o\|^2 + \|S_o'\|^2)\|q\|^2,
$$
$$
\int_DJ(S_{\e}^*, \de\psi )\D\de S \le \frac12\|\D\de S\|^2 + C\lambda_1(\|F\|^2 + \|F'\|^2)\|q\|^2.
$$
Now by {\bf Lemma \ref{lemma}}, we get
$$
\int_D (J(q^1, \psi^1) -  J(q^2, \psi^2))\D\de q \le \|\D\de q\|\|J(q^1, \psi^1) -  J(q^2, \psi^2)\| $$
$$\le \frac12\|\D\de q\|^2 + C(1 + \lambda_1)^2\|\nabla\de q\|^2,
$$
$$
\int_D(J(T^1, \psi^1) -  J(T^2, \psi^2))\D\de T \le \|\D\de T\|\|J(T^1, \psi^1) -  J(T^2, \psi^2)\|$$
$$ \le \frac14\|\D\de T\|^2 + C(\|\nabla\de q\|^2 + |\nabla\de T\|^2),
$$
$$
\int_D(J((1-e^{1-z})\eta_{\e}(z)\theta^1, \psi^1) -  J((1-e^{1-z})\eta_{\e}(z)\theta^2, \psi^2))\D\de T
$$
$$
\le \frac14\|\D\de T\|^2 + C(\|\nabla\de q\|^2 + \|\nabla\de\theta\|^2),
$$
$$
 \int_D(J(S^1, \psi^1) - J(S^2, \psi^2))\D\de S \le \frac12\|\D\de S\|^2 + C(\|\nabla\de q\|^2 + |\nabla\de S\|^2).
 $$
Putting the above estimtates together, we conclude
$$
\frac12\frac{d}{dt}(\|\de\theta_y \|^2 + \|\nabla\de q\|^2 + \|\de T(y, 1)\|^2 + \|\nabla\de T\|^2 + \|\nabla\de S\|^2)
$$
\be\label{cdtotal1}
\le C(\|\de\theta\|^2 + \|\de T(y, 1)\|^2 + \|\nabla\de q\|^2 + \|\nabla\de T\|^2 + \|\nabla\de S\|^2).
\ee
Taking $N$ large enough, multiplying (\ref{ctotal1}) by $N$ and
then adding to (\ref{cdtotal1}), we imply that there exists a
positive constant $\alpha_2$ such that
$$
\frac12\frac{d}{dt}(N\|\de q\|^2 + 2NPrRa^2\lambda_1(\|\de\theta\|^2 + \|\de T\|^2 + \|\de S\|^2))
$$
$$
+ \frac12\frac{d}{dt}(\|\de\theta_y \|^2 + \|\nabla\de q\|^2 + \|\de T(y, 1)\|^2 + \|\nabla\de T\|^2 + \|\nabla\de S\|^2)
$$
$$
\le - \alpha_2(N\|\de q\|^2 + 2NPrRa^2\lambda_1(\|\de\theta\|^2 + \|\de T\|^2 + \|\de S\|^2)) $$
 \be\label{cdtotal2}
 - \alpha_2(\|\de\theta_y\|^2 + \|\nabla\de q\|^2  + \|\nabla\de T\|^2 + \frac12\|\nabla\de S\|^2 + \|\de T(y, 1)\|^2).
\ee
By the Gronwall inequality, we have
$$N\|\de q\|^2 + 2NPrRa^2\lambda_1(\|\de\theta\|^2 + \|\de T\|^2 + \|\de S\|^2)$$
$$ + \|\de\theta_y \|^2 + \|\nabla\de q\|^2 + \|\de T(y, 1)\|^2 + \|\nabla\de T\|^2 + \|\nabla\de S\|^2
$$
$$
\le e^{-\alpha_2 t}(N\|\de q_0\|^2 + 2NPrRa^2\lambda_1(\|\de\theta_0\|^2 + \|\de T_0\|^2 + \|\de S_0\|^2)
$$
\be\label{sctotal}
e^{-\alpha_2 t}(\|\de\theta_{0y}\|^2 + \|\nabla\de q_0\|^2 + \|\de T_0(y, 1)\|^2 + \|\nabla\de T_0\|^2 + \|\nabla\de S_0\|^2).
\ee

This estimate tells us that   any two solution trajectories inside
the absorbing set approach each other as time goes on. This is the
so called strong contraction property.

\begin{remark}
In fact, here we get the estimate is for $\|\nabla\widehat{\de T}\|^2$
and $\|\nabla\widehat{\de S}\|^2 = \|\nabla\de S\|^2$.
But $\|\nabla\widehat{\de T}\|^2 = \|\nabla(\de T + (1 - e^{1-z})\eta_{\e}(z)\de\theta)\|^2$,
$\|\nabla\widehat{\de S}\|^2 = \|\nabla\de S\|^2$ and we have the estimate for
$\|\de\theta\|_1$. Hence we can obtain the similar estimate for
the original $\delta T$.
\end{remark}

Therefore we have the following theorem.

\begin{theorem} \label{contraction}
({\bf Strong contraction property}) Assume that the freshwater
flux $F(y)$ has zero mean as in  (\ref{mean}) and the ocean
basin's latitudinal fraction function $\gamma(y)$ is bounded as in
(\ref{gamma}). Let $a^2 + \|S_a\|^2 + \|S_o\|_{H^2}^2 + \|f\|^2 +
\|F\|_{H^2}^2$ be small enough. Then the coupled atmosphere-ocean
system   (\ref{eqn1})--(\ref{eqn4}) has the strong contraction
property.
\end{theorem}

\bigskip
 \bigskip
Now we come back to the issue of the time-almost periodic (in particular,
time-periodic and time-quasi-periodic)
 motion in the coupled
atmosphere-ocean system.
First, we give the definitions about almost periodic function and pullback
attractor.

A function $\varphi$  $:$ $\rit$ $\to$ $X$, where $(X,d_X)$ is  a
metric space, is called {\em almost periodic\/} \cite{besi} and \cite{sell}
 if for every $\varepsilon$ $>$ $0$  there exists a
relatively dense subset $M_{\varepsilon}$ of $\rit$ such that
$$
d_X \left(\varphi (t+ \tau ), \varphi (t) \right) <  \varepsilon
$$
for all $t$ $\in$ $\rit$ and $\tau$ $\in M_{\varepsilon }$.  A subset $M$
$\subseteq$ $\rit$ is called {\em relatively dense} in $\rit$ if there exists
a positive number $l$ $\in$ $\rit$ such that for every $a$ $\in$ $\rit$ the
interval $[a,a+l]\bigcap \rit$ of length $l$ contains an element of $M$, i.e.
$M\bigcap [a,a+l]$ $\ne$ $\emptyset$ for every $a$ $\in$ $\rit$.

In order to study   the temporally almost periodic solutions of (\ref{eqn1})--(\ref{side}), we need some results from the theory of nonautonomous dynamical systems.
Consider first an autonomous dynamical system on a metric space $P$
described by
a group $\Theta$ $=$ $\{\theta_t\}_{t \in \rit}$ of mappings of $P$ into
itself.

Let $X$ be a complete metric space and consider a continuous mapping
$$
 \Phi : \rit^{+} \times P \times X \to   X
$$
satisfying the properties
$$
\Phi(0,p,\cdot) = {\rm id}_X, \qquad
\Phi(\tau +t,p,x)  =  \Phi(\tau,\theta_t  p, \Phi(t,p,x))
$$
for all $t$, $\tau$ $\in$ $\rit^{+}$, $p$ $\in$ $P$ and $x$ $\in$ $X$.
The mapping $\Phi$ is called a cocycle on $X$ with respect to $\Theta$ on $P$.

The appropriate concept of an attractor for a nonautonomous cocycle systems is
the {\em pullback attractor\/}. In contrast to autonomous attractors it
consists of a family subsets of the original state space $X$ that are
indexed by the cocycle parameter set.

 A family $\widehat{A}$ $=$ $\{A_p\}_{p \in P}$ of nonempty compact sets of
$X$ is
called a {\rm pullback attractor\/} of the cocycle $\Phi$ on $X$ with
respect to
$\theta_t$ on $P$ if it is ${\Phi}$--invariant, i.e.
$$
\Phi(t,p,A_p) = A_{\theta_t} p  \qquad \mbox{for all} \quad t \in \rit^{+},
p \in P,
$$
and {\rm pullback attracting}, i.e.
$$
\lim_{t \to \infty} H^{*}_X\left(\Phi(t,\theta_{-t}p,D), A_p\right) = 0
\qquad \mbox{for all} \quad D \in K(X), \  p \in P,
$$
where $K(X)$ is the space of all nonempty compact subsets of the metric
space $(X,d_X)$. Here $H^{*}_X$ is the Hausdorff
semi--metric between nonempty compact subsets
of $X$, i.e. $H^{*}_X(A,B)$ $:=$ $\max_{a \in A} {\rm dist}(a,B)$ $=$
$\max_{a\in A} \min_{b\in B} d_X(A,b)$ for $A$, $B$ $\in$ $K(X)$.

The following theorem combines  several known results.  See Crauel and
Flandoli
\cite{crflan}, Flandoli and Schmalfu{\ss} \cite{flansch}, and Cheban
\cite{cheban} as well as
\cite{klosch, cdf, chklsch}  for this and various related proofs.

\begin{theorem} \label{th1}
Let $\Phi$ be a continuous cocycle on a metric space $X$ with respect to a
group
$\Theta$ of continuous mappings on a metric space $P$. In addition, suppose
that
there is a nonempty compact subset $B$ of $X$ and that for
every $D$ $\in$ $K(X)$
there exists a $T(D)$ $\in$ $\rit^{+}$, which is independent of $p$ $\in$ $P$,
 such that
\begin{equation}\label{fa}
\Phi(t,p,D) \subset B \quad \mbox{for all} \quad t > T(D).
\end{equation}
Then  there exists a unique pullback attractor
$\widehat{A}$ $=$
$\{A_p\}_{p \in P}$ of the cocycle $\Phi$ on $X$, where
\begin{equation}\label{pbat}
A_p = \bigcap_{\tau \in \rit^{+}} \overline{\bigcup_{t > \tau \atop t \in
\rit^{+}}
\Phi\left(t,\theta_{-t}p,B\right)}.
\end{equation}
Moreover, if the cocycle $\Phi$  is strongly
contracting inside
the absorbing set $B$. Then the pullback attractor consists of  singleton
valued sets,
i.e.   $A_p$ $=$ $\{a^*(p)\}$,  and the
mapping $p$ $\mapsto$ $a^*(p)$ is continuous.
\end{theorem}

The  solution operators $S_{t,t_0}$ for (\ref{eqn1})--(\ref{side}) form a cocycle
mapping on $X = H^1(0, 1)\times H_0^1(D)\times H^1(D) \times H^1(D)$ with
parameter set $P$ $=$ $\rit$, where $p$ $= t_0$, the initial time, and
$\theta_t t_0$ $=$ $t_0+t$, the left shift by time $t$. However, the space
$P$ $=$
$\rit$ is not compact here. Though more complicated, it is more useful to
consider
$P$ to be the closure of the subset $\{\theta_t f, t \in \rit\}$, i.e. the
hull of $f$,  in
the metric  space $L^2_{loc}\left(\rit,X)\right)$ of locally
$L^2(\rit)$--functions
$f$ $:$ $\rit$ $\to$ $X$ with the metric
$$
d_P(f,g) :=  \sum_{N=1}^{\infty} 2^{-N} \min\left\{1, \sqrt{\int_{-N}^N
|\|f(t)-g(t)\||^2 \, dt} \right\}
$$
with  $\theta_t$ defined to be  the left shift operator, i.e.
$\theta_t f(\cdot) := f(\cdot+t)$, where $|\|\cdot\||$ denotes the norm in $X$.
By a classical result \cite{besi,sell},  a function $f$ in the
above metric space is almost periodic if and only if the the  hull of $f$
is compact
and minimal. Here minimal means nonempty, closed and invariant with respect
to the
autonomous dynamical system generated by the shift operators $\theta_t$
such that
with no proper subset has these properties. The cocycle mapping is defined
to be the
solution $\vec{\omega}(t) = \{\theta, q, T, S\}$ of (\ref{eqn1})--(\ref{side}) starting at $\vec{\omega}_0 = \{\theta_0, q_0, T_0, S_0\} \in X$ at time $t_0$ $=$ $0$ for a given
forcing mapping $f$ $\in$ $P$, i.e.
$$
\Phi(t,f,\om_0) := S_{t,0}^f \ \om_0,
$$
where we have included a superscript $f$ on $S$ to denote the dependence on
the forcing
term $f$. (This dependence is in fact continuous, see \S 3). The cocycle property
here follows from the fact  that
$S_{t,t_0}^f \vec{\om}_0$ $=$ $S_{t-t_0,0}^{\theta_{t_0}f}\ \vec{\om}_0$ for
all $t$ $\geq$ $t_0$, $t_0$ $\in$ $\rit$,
$\vec{\om}_0 \in X$
and $f$ $\in$ $P$.

\bigskip

Following   Theorem 5.3 and
the dissipativity and
strong contraction properties shown in the  last two sections,
we conclude that  the coupled atmosphere-ocean system
 (\ref{eqn1})--(\ref{eqn4}) has
a unique pullback attractor,
  consists of the  singleton valued component
  $\{\vec{a}^*(p)\} \in \hat{A}$ and the mapping
$p$ $\mapsto$ $\vec{a}^{*}(p)$ is continuous on $P$.
As in Duan and Kloeden \cite{duklo} or Gao, Duan and Fu \cite{gdf},
we now show that this
singleton attractor $\vec{a}^{*}(p)$ defines an almost periodic
solution.

In fact, the mapping  $p$ $\mapsto$ $\vec{a}^{*}(p)$
is uniformly continuous on $P$ because $P$ is compact subset of
$L^2_{loc}\left(\rit, X)\right)$ due to the assumed  almost periodicity.
That is, for every $\varepsilon$ $>$ $0$ there exists a $\delta(\varepsilon)$
$>$ $0$ such that $\|\vec{a}^*(p) - \vec{a}^*(q) \|$ $<$ $\varepsilon$ whenever
$d_P(p,q)$ $<$
$\delta$.  Now let  the point $\bar{p}$ ($=$ $f$, the given temporal
forcing function)
be almost periodic and   for $\delta$ $=$ $\delta(\varepsilon)$
$>$ $0$ denote by $M_{\delta}$ the relatively
dense subset of $\rit$ such that
$d_P(\theta_{t+\tau}\bar{p},\theta_t \bar{p})$ $<$ $\delta$
for all $\tau$ $\in$ $M_{\delta}$ and
$t$ $\in$ $\rit$. From this and the uniform
continuity we have
$$
\|\vec{a}^*(\theta_{t+\tau} \bar{p}) - \vec{a}^*(\theta_t \bar{p})\| < \varepsilon
$$
for all $t$ $\in$ $\rit$ and $\tau$ $\in$ $M_{\delta(\varepsilon)}$. Hence
$t$ $\mapsto$ $\{\theta^*, q^*, T^*, S^*\}(t)$ $:=$  $\vec{a}^*(\theta_{t}\bar{p})$ is almost periodic, and
it is a solution of the  coupled atmosphere-ocean model.
It is unique as the single-trajectory pullback attractor is the only trajectory
that exists and is bounded for the entire time line.
Therefore, we have the following result.

\begin{theorem} \label{periodic}
({\bf Periodic, quasiperiodic and almost periodic motion})
Assume that the freshwater flux $F(y)$ has zero mean as in  (\ref{mean})
and the ocean basin's latitudinal fraction function $\gamma(y)$
is bounded as in (\ref{gamma}).
Let $a^2 + \|S_a\|^2 + \|S_o\|_{H^2}^2
+ \|f\|^2 + \|F\|_{H^2}^2$ be small enough.
Then the coupled atmosphere-ocean system   (\ref{eqn1})--(\ref{eqn4})
has  unique time-periodic, quasiperiodic and almost periodic motions,
 when the external fluctuation $f$ in the atmospheric energy balance model
  is time-periodic, quasiperiodic and almost periodic, respectively.
\end{theorem}

This result   may be relevant to the El Nino-Southern Oscillation
phenomenon.
Ei Nino is a well-known climate phenomenon in the  atmosphere-ocean
system. Originally, it refers to a seasonal invasion, along the coast of Peru around Christmas, of a warm southward ocean current that displaced the north-flowing cold current. Today, it is regarded as a
part of a  phenomenon called El Nino-Southern Oscillation (ENSO),
a continual but  ``quasi"-periodic, perhaps irregular,
cycle of shifts in ocean and atmosphere condition that affect the globe
climate \cite{Neelin, DijkBook, Siedler}.

%%%%%%%%%%%%%%%%%%%%%%
\section{Summary}

We have investigated the dynamical behavior of a
 coupled atmosphere-ocean system.

First, we show that the coupled atmosphere-ocean system is stable
under the external fluctuation in the atmospheric energy balance relation
(Theorem \ref{stability}).
Then, we estimate the
atmospheric temperature feedback in terms of the freshwater flux, heat flux
and the   external fluctuation at the air-sea interface, as well as
the earth's longwave radiation
and the shortwave solar radiation (Theorem \ref{feedbacktheorem}).
Finally, we prove that the coupled atmosphere-ocean system
has    time-periodic, quasiperiodic and almost periodic motions
(under suitable conditions on the physical quantities such as
freshwater flux, the earth's longwave  radiative
cooling coefficient and the shortwave
solar   radiation profile),
when the external fluctuation  in the atmospheric energy balance relation
is time-periodic, quasiperiodic and almost periodic, respectively
(Theorem \ref{periodic}).

%%%%%%     Acknowledgement          %%%%%
%%%%%%%%%%%%%%%%%%%%%%%%%
%%%%%%%%%%%%%%%%%%%
\bigskip

{\bf Acknowledgement.}  A part of this work was done at the
Oberwolfach Mathematical Research Institute, Germany. This work
was partly supported by the NSF Grant DMS-0209326,   the Grant
10001018 of the NNSF of China, and the Grant BK2001108 of the  NSF
of Jiangsu Province, as well as the Scientific Research Foundation
for Returned Overseas Chinese Scholars of Jiangsu Education
Commission. And a part of this work was done while H. Gao was
visiting the  Illinois Institute of Technology and the Institute
for Mathematics and Its Applications, USA. This research was
supported in part
       by the Institute for Mathematics and its Applications with funds
       provided by the National Science Foundation.

\end{document}